\documentclass[final,12pt]{article}
\usepackage{amsfonts,amssymb,amsmath,amsthm,graphicx,overpic,color}
\usepackage{xcolor,soul,cancel}
\usepackage[color]{showkeys}
\newtheorem{thm}{Theorem}

\newtheorem{lem}{Lemma}
\newtheorem{pro}{Proposition}

\newcommand{\conv}{\operatorname{conv}}

\newcommand{\rk}{\operatorname{rank}}
\newcommand{\scal}[1]{\langle#1\rangle}

\newcommand{\cn}{\operatorname{cn}}

\newcommand{\sn}{\operatorname{sn}}
\newcommand{\eq} [1] {\begin{equation}\label{#1}\quad}
\newcommand{\en} {\end{equation}}

\newcommand{\FigDir}{}

\definecolor{darkgreen}{rgb}{0,0.5,0}
\definecolor{darkred}{rgb}{0.7,0,0}

\title{On partial isometries with circular numerical range}

\author{Elias Wegert and Ilya Spitkovsky}

\begin{document}
\maketitle

\begin{abstract}
In their LAMA'2016 paper Gau, Wang and Wu conjectured that a partial isometry
$A$  acting on $\mathbb{C}^n$ cannot have a circular numerical
range with a non-zero center, and proved this conjecture for
$n\leq 4$. We prove it for operators with $\rk A=n-1$ and
any $n$.

The proof is based on the unitary similarity of $A$ to a compressed
shift operator $S_B$ generated by a finite Blaschke product $B$.
We then use the description of the numerical range of $S_B$ as
intersection of Poncelet polygons, a special representation of
Blaschke products related to boundary interpolation, and an explicit
formula for the barycenter of the vertices of Poncelet polygons involving
elliptic functions.
\end{abstract}

\section{Introduction} \label{sec.Intro}
Denote by $\mathbb{C}^{n\times n}$ the algebra of all $n$-by-$n$ matrices with
complex entries. The numerical range $W(A)$ of $A\in\mathbb{C}^{n\times n}$ is
the set of values of the quadratic form $\scal{Ax,x}$ on the unit sphere of
$\mathbb{C}^n$. By the celebrated Toeplitz-Hausdorff theorem, $W(A)$ is a
convex subset of $\mathbb{C}$; see \cite{GusRa}, Chapter~1 of \cite{HJ2}, or
the recent books~\cite{DGSV},~\cite{GaMaRo} for this and other properties
of the numerical range.

It is easy to see that  $W(A)$ contains the spectrum $\sigma(A)$ of $A$, and
therefore its convex hull $\conv\sigma(A)$. The two sets coincide for normal
matrices $A$, but not in general. So, for a unitary matrix $U$ the set $W(U)$
is a polygon inscribed into the unit circle.

In this paper we are concerned with {\em partial isometries},
i.e. matrices $A$ the action of which preserves norms of vectors from $(\ker
A)^\perp$. Every partial isometry $A$ is the orthogonal sum of unitarily
irreducible partial isometries, a unitary component $U$, and a zero block (with
each of the components allowed to be missing).
Since the numerical range of a block diagonal matrix is the convex hull of the
numerical ranges of its blocks, the emphasis can be put on consideration of
unitarily irreducible partial isometries.

A canonical example of the latter is the Jordan block $J_n$. Observe that
$W(J_n)$ is the circular disk centered at the origin of the radius
$\cos\frac{\pi}{n+1}$. Simple examples show that, even for $n=2$, the numerical
range of a partial isometry is not necessarily a circular disk. However, for
$n\leq 4$ it was observed by Gau, Wang and Wu in \cite{GWW16} that \emph{if}\,
$W(A)$ \emph{happens to be a circular disk}, then this disk is necessarily centered
at the origin. The authors conjectured that this property persists for all $n$.

The case $n=5$ was settled in our recent (joint with I. Suleiman) paper~\cite{SSW}
by a straightforward but rather lengthy proof based on divisibility
considerations of the so-called Kippenhahn polynomial.
It was also observed there that it follows from the results
of~\cite{GWW16} that the conjecture holds for $A$ of rank one or two,
independent of the value of $n$.

Another extreme situation, when $A$ has a one-dimensional kernel, is
non-trivial. The respective conjecture was stated separately in~\cite{GWW16};
for convenience of reference we will call it the {\em special}
Gau-Wang-Wu conjecture. To prove it for all values of $n$ is the goal of this
paper.

\begin{thm}[]\label{thm.GWW}
Let $A$ be a partial isometry acting on $\mathbb{C}^n$. If $\dim\ker A=1$ and
the numerical range $W(A)$ of $A$ is a circular disk, then this disk is centered
at the origin.
\end{thm}

Note that a unitarily reducible partial isometry $A$ with one-dimensional kernel
is the orthogonal sum of a unitarily irreducible one and a non-trivial unitary
matrix. As such, $W(A)$ cannot be a circular disk. So, in the rest of the paper
we concentrate on unitarily irreducible matrices only.


\section{Blaschke, Kippenhahn and Poncelet} \label{sec.KBP}

In this section we summarize some relevant facts from operator theory,
complex analysis, and geometry. For more detailed information
we refer to~\cite{GWW16} and the books~\cite{DGSV},\cite{GaMaRo}.


Recall first that Kippenhahn \cite{Kipp1951} describes the numerical range
of operators $A$ acting in $\mathbb{C}^n$ as the convex hull of (the real part
of) an algebraic curve $C(A)$ of class $n$, often called the \emph{Kippenhahn curve}
of $A$ (see~\cite{Ki08} for an English translation, and \cite[Section~13]{DGSV}
for a contemporary treatment).


An operator $A$ acting in $\mathbb{C}^n$ is a unitarily irreducible partial
isometry with $\dim\ker A=1$
if and only if $A$ is a non-invertible matrix of class $\mathcal{S}_n$,
which consists of the contractions $A\in\mathbb{C}^{n\times n}$ with
eigenvalues in the unit disk $\mathbb{D}$ and $\rk (I-A^*A)=1$
(see \cite[Proposition 2.3]{GWW16}).
Specific properties of the numerical range of operators in $S_n$
were studied independently by Gau and Wu~\cite{GauWu98} and
Mirman~\cite{Mir} (see \cite{DaGoVo} for further information).
Gau and Wu prove that among all operators acting in $\mathbb{C}^n$ those
in $S_n$ are distinguished by the so-called \emph{Poncelet property} of
their numerical range: For each $t$ on the unit circle $\mathbb{T}$
there exists a $(n+1)$-gon $P_t$ which is inscribed in $\mathbb{T}$,
circumscribed about $W(A)$, and has $t$ as a vertex. The vertices of
these \emph{Poncelet polygons} $P_t$ are the eigenvalues
of unitary dilations $U_t$ of $A$ (\cite[Theorem~2.1]{GauWu98},
for an alternative proof see~\cite{DaGoVo}).


Operators in $\mathcal{S}_n$ have simple models which we describe next.
To begin with, let $B$ be a Blaschke product of degree $n$ with zeros
$\lambda_1,\ldots,\lambda_n$ in $\mathbb{D}$,
\begin{equation} \label{eq.DefBz}
B(z):=\gamma\,\prod_{k=1}^{n}\frac{z-\lambda_k}{1-\overline{\lambda_k}z},
\qquad |\gamma| = 1.
\end{equation}
What follows is independent of the unimodular factor $\gamma$,
so that we often assume $\gamma=1$.
The {\em model space} ${\mathcal K}_B$ of $B$ is the $n$-dimensional
linear space of all rational functions $p/q$ with denominator
$q(z):=\prod(1-\overline{\lambda_k}z)$ and $\deg p\le n-1$
(see \cite[Chapter~12]{GaMaRo} or \cite{GoWa17}, for instance).
The \emph{compressed shift} $S_B$ generated by $B$ is the operator acting
in ${\mathcal K}_B$ as the compression of the multiplication by $z$ into
${\mathcal K}_B$ (by orthogonal projection of $L^2(\mathbb{T})$ onto its
subspace ${\mathcal K}_B$).

The notation $\lambda_k$ for the zeros of $B$ was chosen intentionally:
$\lambda_1,\ldots, \lambda_n$ are the eigenvalues of $S_B$ (\cite[Corr.12.6.7]{GaMaRo}).
Moreover, the numerical range $W(S_B)$ of $S_B$ has a beautiful
geometric description, which reflects Kippenhahn's theorem as well as the
Poncelet property
(see~\cite{ChaGoPaRo18}, 
\cite[p.288]{GaMaRo},
\cite{GoWa17}). 

\begin{pro}[]\label{pro.NumRanShft}
Let $B$ be a Blaschke product of degree $n$ with zeros $\lambda_1,\ldots,\lambda_n$
and define $B_1$ by $B_1(z):=z \,B(z)$. For $t \in\mathbb{T}$ let $P_t$ be the
convex $(n+1)$-gon with vertices at the preimages $B_1^{-1}(t)$ of $t$.
Then the numerical range of $S_B$ is
\[
W(S_B) = \bigcap_{t\in\mathbb{T}}\ \conv P_t.
\]
\end{pro}

\noindent
The sides of all \emph{Poncelet polygons} $P_t$ are tangent to the \emph{Kippenhahn
curve} of $S_B$; this curve is generated as an envelope of straight lines connecting
successive points on $\mathbb{T}$ at which $B_1$ has constant
phase $B_1/|B_1|=t$. Since this curve is solely defined by the Blaschke product
$B_1$ and has the Poncelet property, we call it the \emph{Poncelet curve}
of $B_1$. The numerical range $W(S_B)$ is the closure of the interior of that curve.

Figure~1 illustrates this construction in the ``phase plots''
of two Blaschke products $B_1(z)=z\,B(z)$ with degree 4 and 5, respectively.
The functions are depicted on their domain $\mathbb{D}$, coloring a point $z$
according to the phase $B_1(z)/|B_1(z)|$ of the function value.
The points where all colors meet are the zeros of $B_1$. Those zeros different
from $0$ are the eigenvalues of $S_B$. For more detailed explanations
of phase plots we refer to~\cite{WeSe} and~\cite{Wgt2012}.

\begin{center}
\includegraphics[width=0.45\linewidth]{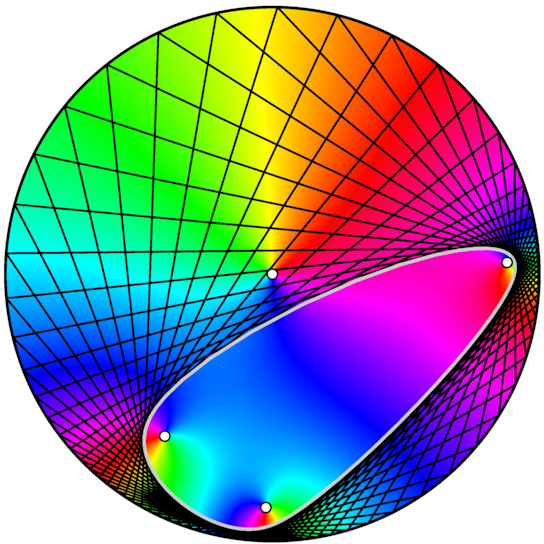}
\hspace{0.05\linewidth}
\includegraphics[width=0.45\linewidth]{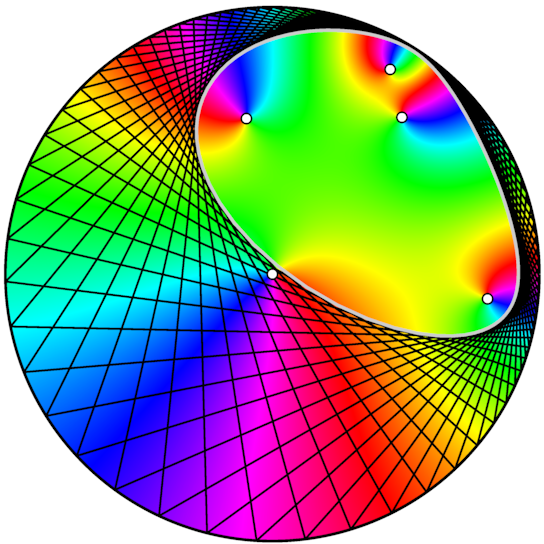}\\[1ex]
Figure~1. Generation of the Poncelet curve for Blaschke products
of degree $4$ and $5$.
\end{center}


\noindent
Compressed shift operators are the typical representatives of
the class $S_n$: If $A\in S_n$ has the eigenvalues $\lambda_1,\ldots,\lambda_{n-1},
\lambda_n=0$, and $B$ is its \emph{associated Blaschke product} defined by
\begin{equation} \label{eq.DefBz2}
B(z):=z\,\prod_{k=1}^{n-1}\frac{z-\lambda_k}{1-\overline{\lambda_k}z},
\end{equation}
then $A$ is unitarily similar to $S_B$
(see, e.g.,~\cite[Theorem 12.7.8]{GaMaRo}). So it suffices to verify the first
Gau-Wang-Wu conjecture for the operators $S_B$ generated by Blaschke products
\eqref{eq.DefBz2}.
Since these Blaschke products have a zero at $0$, the Blaschke products
$B_1(z):=z\,B(z)$ have a \emph{double zero} at the origin.
What remains to prove is that the Poncelet curve $C(B_1)$ can only be circular
if it is centered at the origin.


\section{Boundary representation of Blaschke products} \label{sec.abc}

A crucial ingredient to the proof of Theorem~\ref{thm.GWW} is a special
representation of Blaschke products. The usual way of writing these functions
as in \eqref{eq.DefBz}
emphasizes the role of their zeros.  Since Blaschke products are objects
of hyperbolic geometry (often called ``hyperbolic polynomials''),
and the origin is not a distinguished point in that geometry,
one may ask for alternative descriptions. In this section we propose
a representation that uses values of $B$ on the unit circle
and is related to special \emph{boundary interpolation problems}.

To begin with, we rewrite $B$ from \eqref{eq.DefBz} as rational function $B=p/q$,
with polynomials
\[
p(z):= d \,\prod_{k=1}^{n} (z-\lambda_k)
= \sum\limits_{k=0}^{n} p_k z^k, \qquad
q(z):= z^n \overline{p(1/\overline{z})}
= \sum\limits_{k=0}^{n} \overline{p}_{n-k} z^k,
\]
and $d^2=\gamma$ (it does not matter which square root of $\gamma$ we chose).
Then we pick three pairwise distinct unimodular complex numbers $a,b,c$
which we assume to be cyclically ordered on the unit circle $\mathbb{T}$
such that $a \prec b \prec c \prec a$
(for example $a=1$, $b=-1$, $c=-\mathrm{i}$). The sets $A_n:=B^{-1}(a)$ and
$B_n:=B^{-1}(b)$ of pre-images of $a$ and $b$ consist of $n$ pairwise distinct
points $a_1,\ldots,a_n$ and $b_1,\ldots,b_n$, respectively.
Since the argument of $B(\mathrm{e}^{\mathrm{i}\theta})$ is a strictly increasing
function of $\theta$, the points in $A_n$ and $B_n$ must be interlacing on
$\mathbb{T}$, so that we may assume the cyclic ordering
\begin{equation} \label{eq.4}
a_1 \prec b_1 \prec a_2 \prec b_2 \prec \ldots a_n \prec b_n \prec a_{n+1}:=a_1.
\end{equation}
Moreover, each positively oriented arc $(b_k,a_{k+1})$ from $b_k$ to $a_{k+1}$
contains exactly one point $c_k$ with $B(c_k)=c$. We choose just one of them
and denote it by $c_0$.
So, given $a,b,c \in \mathbb{T}$, the Blaschke product $B$ defines $2n+1$
points $a_1,\ldots, a_n$, $b_1, \ldots, b_n$ and $c_0$ on $\mathbb{T}$.

This construction also works the other way around. The determination of $B$
from the sets $A_n$, $B_n$ and the point $c_0$ requires the solution
of the interpolation problem
\begin{equation} \label{eq.5}
B(a_k)=a,\quad B(b_k)= b, \quad B(c_0)=c, \qquad k=1,\ldots,n.
\end{equation}
While \emph{Nevanlinna-Pick} interpolation problems $B(z_k)=w_k$ with
$|z_k|,|w_k| <1$ are studied and understood for more than a century,
the history of interpolation problems with $|z_k|=|w_k|=1$
is much shorter (two milestones are Cantor and Phelps~\cite{CaPh1965},
Jones and Ruscheweyh~\cite{JoRu1987}).
Since then these problems have attracted quite some interest, but a number
of questions is still unanswered. In particular, no algebraic criterion for
the determination of the minimal degree of an interpolant seems to be known
(see Semmler and Wegert~\cite{SeWe2006}, Glader~\cite{Gla2008},
and Bolotnikov~\cite{Bol2018}). Fortunately, the problem at hand is fairly
well understood.

\begin{thm}[]\label{pro.1}
Assume that $a_1,\ldots,a_n$ and $b_1,\ldots,b_n=:b_0$ are points on the unit
circle $\mathbb{T}$, strictly cyclically ordered according to \eqref{eq.4}, and
let $c_0\in \mathbb{T}$  be such that $b_{k-1}\prec c_0 \prec
a_k$ for some $k$. Then, for any triple $a,b,c$ of points on $\mathbb{T}$,
cyclically ordered such that $a \prec b \prec c \prec a$, there exists a unique
Blaschke product $B$ of degree $n$ that satisfies \eqref{eq.5}.
\end{thm}

\begin{proof}
Corollary~10 in Daepp, Gorkin and Voss \cite{DaGoVo}  tells us that
there exists a Blaschke product $B_0$ of degree $n$ such that $B_0(a_k)=a_0$
and $B_0(b_k)=b_0$ for some $a_0,b_0 \in\mathbb{T}$ and $k=1, \ldots, n$.
The assumptions on the ordering of $a_k$, $b_k$ and $c_0$ (and the fact that
$B_0$ is an orientation preserving $n$-fold covering map of $\mathbb{T}$ onto
itself), guarantee that the triple $\big(a_0,b_0,B(c_0)\big)$ has the same
orientation as $(a,b,c)$. The composition $B:=B_1 \circ B_0$ with the Blaschke
factor $B_1$ that maps $a_0\mapsto a$, $b_0 \mapsto b$ and $B(c_0)\mapsto c$
is a solution of the interpolation problem.

Since $B$ has degree $n$ and satisfies $2n+1$ interpolation conditions,
the interpolation problem falls in the class of ``elastic'' problems,
and uniqueness follows from~\cite[Theorem~1]{SeWe2006}.
\end{proof}

\noindent
In order to construct the solution explicitly, we define the polynomials
\begin{equation} \label{eq.6}
P := q - \overline{a}p, \quad Q:= q - \overline{b}p,
\end{equation}
so that
\begin{equation} \label{eq.7}
B = \frac{p}{q} = \frac{P-Q}{\overline{b}P-\overline{a}Q}.
\end{equation}
Clearly we have
\begin{equation} \label{BBB}
B=a\ \Leftrightarrow\ P=0,\qquad
B=b\ \Leftrightarrow\ Q=0,\qquad
B=0\ \Leftrightarrow\ P=Q.
\end{equation}
Hence, setting
\begin{equation} \label{eq.PQ}
\widetilde{P}(z) := \prod_{k=1}^{n}(z-a_k), \quad
\widetilde{Q}(z) := \prod_{k=1}^{n}(z-b_k),
\end{equation}
we conclude that $P=\alpha\,\widetilde{P}$ and $Q = \beta\,\widetilde{Q}$
with some non-zero numbers $\alpha$ and $\beta$. From $B(c_0)=c$ it follows that
\[
\alpha\,\widetilde{P}(c_0)- \beta\,\widetilde{Q}(c_0)
= \alpha\,\overline{b}c\,\widetilde{P}(c_0) - \beta\,\overline{a}c\,\widetilde{Q}(c_0),
\]
which is satisfied for
\begin{equation} \label{eq.9}
\alpha := (1-\overline{a}c)\,\widetilde{Q}(c_0), \quad
\beta  := (1-\overline{b}c)\,\widetilde{P}(c_0).
\end{equation}
Note that the numbers $\alpha$ and $\beta$ are uniquely determined up to a common
factor. Then we have
\begin{equation} \label{eq.10}
\widetilde{P}(c_0) = \prod_{k=1}^{n}(c_0-a_k), \quad
\widetilde{Q}(c_0) = \prod_{k=1}^{n}(c_0-b_k),
\end{equation}
and the explicit solution of the boundary interpolation problem~\eqref{eq.5}
is given by
\begin{equation} \label{eq.11}
B = \frac{\alpha \widetilde{P}-\beta \widetilde{Q}}
{\overline{b} \,\alpha \widetilde{P}-\overline{a}\,\beta \widetilde{Q}},
\end{equation}
with $\widetilde{P}$, $\widetilde{Q}$ from \eqref{eq.PQ}, and $\alpha$, $\beta$
from \eqref{eq.9} and  \eqref{eq.10}.
This is the desired boundary representation of $B$.


\section{Proof of the special Gau-Wang-Wu conjecture}
\label{sec.Main}

With the help of Proposition~\ref{pro.NumRanShft}, Theorem~\ref{thm.GWW} can
be recast as follows. Note that, in this section, $B$ stands for the Blaschke product
formerly denoted by $B_1$ and $n:=\deg B$.

\begin{thm}[]\label{thm1}
Let $B$ be a Blaschke product of degree $n\ge 3$ with $B(0)=0$ and $B'(0)=0$.
If the Poncelet curve associated with $B$ is a circle $C$, then its center
$c$ is the origin.
\end{thm}

The proof will occupy the rest of this section.
To begin with, we observe that the center $c$ of the circle $C$
determines the Blaschke product $B$ with $B(0)=0$ almost uniquely.
As we shall show, this Blaschke product satisfies $B'(0)=0$ if and only
if $c=0$. Interestingly, this can be reduced to a problem of plane geometry.

For a fixed center $c\in\mathbb{D}$ of $C$ there is a unique radius $r$ such that
the circle $C$ has a circumscribed $n$-gon with vertices on the unit circle
$\mathbb{T}$. We assume that $c>0$ and fix the corresponding radius $r$.

As we have seen in Section~\ref{sec.KBP}, for each $t\in \mathbb{T}$ the
preimages $B^{-1}(t)$ are the vertices of a Poncelet $n$-gon $P_t$
circumscribed about $C$.
Among all these polygons there are exactly two which are symmetric
with respect to the real line: one with a vertex at $-1$, and a second
one with a side (``on the left'') parallel to the imaginary axis
(see Figure~2).
\bigskip\par

 \begin{center}
 \begin{overpic}[width=.4\linewidth,tics=5]{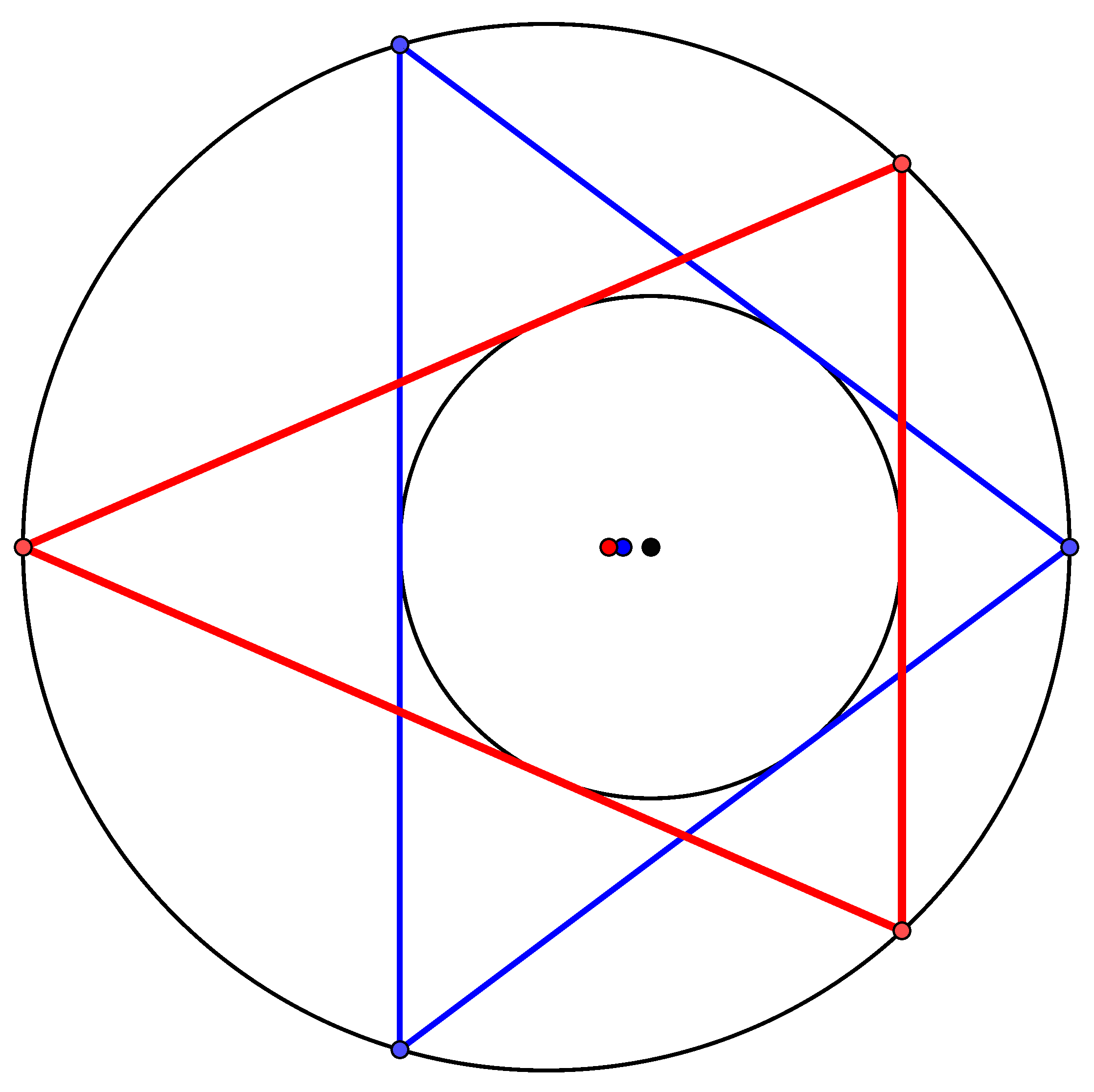}
 \fontsize{12pt}{14pt}\selectfont
 \put(0,50){\makebox(0,0)[rc]{$a_1$}} %
 \put(85,10){\makebox(0,0)[lc]{$a_2$}}
 \put(85,90){\makebox(0,0)[lc]{$a_3$}}
 \put(36,0){\makebox(0,0)[cc]{$b_1$}}
 \put(101,50){\makebox(0,0)[lc]{$b_2$}}
 \put(36,100){\makebox(0,0)[cb]{$b_3$}}
 \put(51,43){\makebox(0,0)[cb]{$a_*$}}
 \put(57,43){\makebox(0,0)[cb]{$b_*$}}
 \put(63,44){\makebox(0,0)[cb]{$c$}}
 \end{overpic}
 \hspace{.1\linewidth}
 \begin{overpic}[width=.4\linewidth,tics=5]{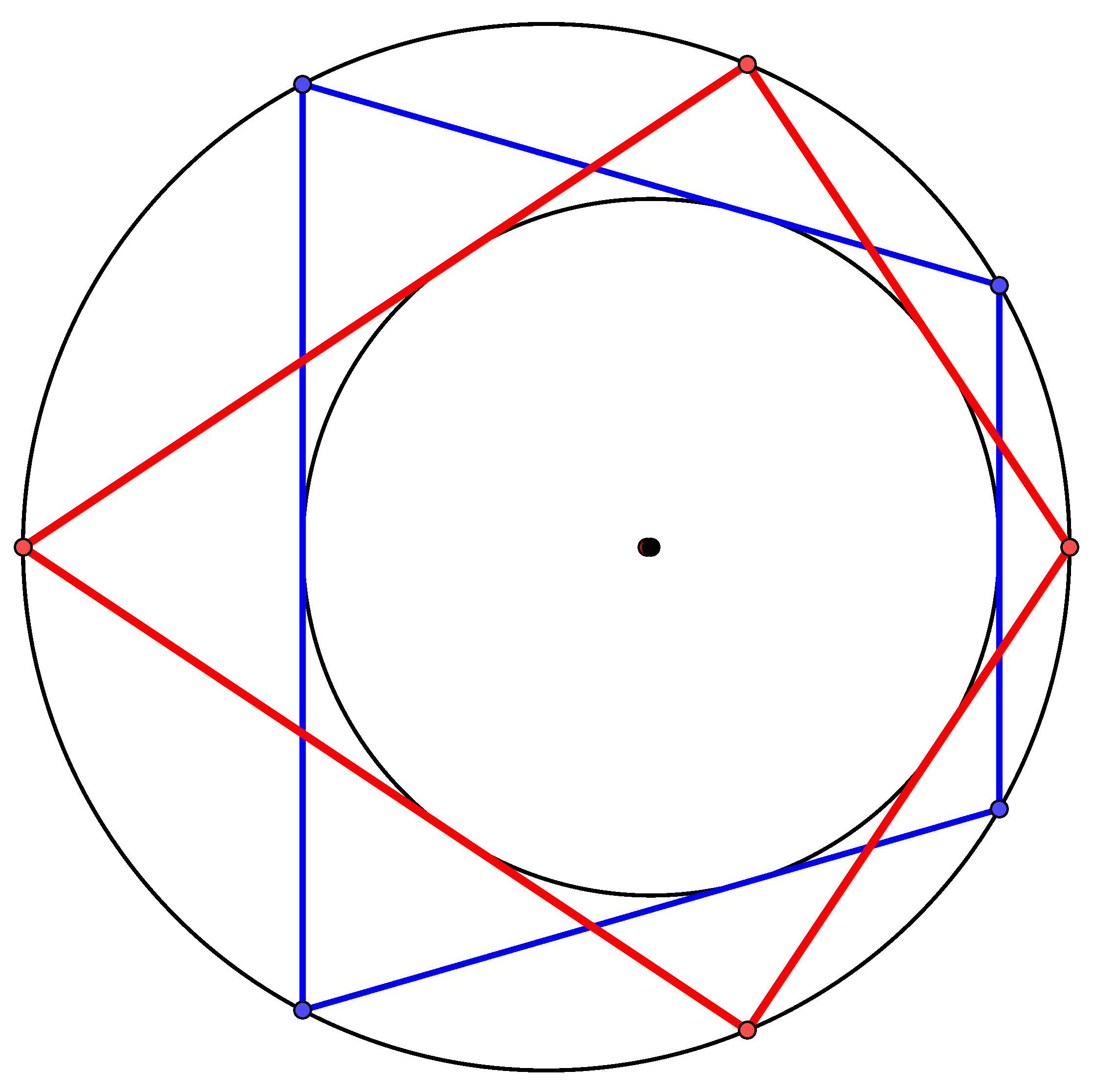}
 \fontsize{12pt}{14pt}\selectfont
 \put(0,50){\makebox(0,0)[rc]{$a_1$}}
 \put(68,0){\makebox(0,0)[lc]{$a_2$}}
 \put(100,50){\makebox(0,0)[lc]{$a_3$}}
 \put(68,100){\makebox(0,0)[lc]{$a_4$}}
 \put(28,4){\makebox(0,0)[ct]{$b_1$}}
 \put(95,22){\makebox(0,0)[lc]{$b_2$}}
 \put(95,78){\makebox(0,0)[lc]{$b_3$}}
 \put(28,96){\makebox(0,0)[cb]{$b_4$}}
 \put(56,50){\makebox(0,0)[rc]{$a_*$}}
 \put(62,50){\makebox(0,0)[lc]{$b_*=c$}}
 \end{overpic}
 \bigskip\par
 Figure~2. The symmetric Poncelet polygons for $c=0.15$, $n=3$ and $n=4$.
 \end{center}
\smallskip\par
\noindent
The vertices $a_k$ and $b_k$ of these polygons form two interlacing sets
on the unit circle $\mathbb{T}$, and we may assume that they are cyclically
ordered,
\[
-1=a_1 \prec b_1 \prec a_2 \prec b_2 \prec \ldots \prec a_n \prec b_n \prec -1.
\]
According to the general property of Poncelet polygons associated with Blaschke
products we have
\begin{equation} \label{eq.NRD1}
B(a_k)=a, \quad B(b_k)=b, \qquad k=1,\ldots n,
\end{equation}
for some $a,b \in\mathbb{T}$ and $a \not=b$. This is exactly the situation we
have encountered in the preceding section. Representing $B$ as
$(P-Q)/(\overline{b}P-\overline{a}Q)$ as in \eqref{eq.7},  from
\eqref{BBB} we have
\begin{equation} \label{eq.8}
P(z) = \alpha\prod_{k=1}^{n}(z-a_k), \quad
Q(z) = \beta\prod_{k=1}^{n}(z-b_k), \qquad
\alpha,\beta \in \mathbb{C}\setminus \{0\}.
\end{equation}
The assumption $B(0)=0$ implies that $P(0)=Q(0)$. Moreover,
$B'=(P'Q-PQ')/Q^2$, so that $B'(0)=0$ if and only if $P'(0)=Q'(0)$.
Since
\[
P'(z)=P(z)\,\sum\limits_{k=1}^{n} \frac{1}{z-a_k}, \qquad
Q'(z)=Q(z)\,\sum\limits_{k=1}^{n} \frac{1}{z-b_k},
\]
and $P(0)=Q(0)\not=0$, this is equivalent to
\[
\sum\limits_{k=1}^{n} \frac{1}{a_k} = \sum\limits_{k=1}^{n} \frac{1}{b_k}.
\]
Using $|a_k|=|b_k|=1$, as well as the symmetries of the vertex sets
with respect to $\mathbb{R}$, we get
\begin{equation} \label{eq.NRD2}
\sum\limits_{k=1}^{n} a_k = \sum\limits_{k=1}^{n} \overline{a}_k
=\sum\limits_{k=1}^{n} \frac{1}{a_k}
=\sum\limits_{k=1}^{n} \frac{1}{b_k}
= \sum\limits_{k=1}^{n} \overline{b}_k = \sum\limits_{k=1}^{n} b_k.
\end{equation}
This equation has a nice geometric interpretation: The barycenters
(centers of mass) $a_*$ and $b_*$ of the vertices $a_k$ and $b_k$ of
the two Poncelet $n$-gons must coincide. We will prove that this can
only happen if $c=0$.

A recent paper by Richard Schwartz and
Sergei Tabachnikov~\cite{SchT} studies the locus of the barycenters of all
Poncelet polygons inscribed in and circumscribed about ellipses. We quote their
main result, adapted to our situation.%
\footnote{The authors attribute this result to some Konstantin Shestakov,
who served for the Russian army in the war against Napoleon. Though this
person and his story are very likely inventions, we warmly recommend to
read this masterpiece of fictitious history.}

\begin{thm}[Schwartz and Tabachnikov]\label{thm.SST}
The locus $S$ of the barycenters of the vertices of the Poncelet polygons
$P_t$ is a circle or a point.
\end{thm}

It follows from symmetry arguments that $S$ is symmetric with respect to the
real line. Assuming that it is not a point, $S \cap \mathbb{R}=\{a_*,b_*\}$
and $a_*\not=b_*$. To exclude that $S$ is a point, one
could analyze the proof given by Schwartz and Tabachnikov. Experts in
projective geometry who see this immediately may skip the rest of the paper.
For those less familiar with these techniques, we adopt Jacobi's traditional
approach to Poncelet's theorem for two
circles using elliptic functions (Jacobi~\cite{Jac1828}, see Griffith~\cite{Gri}
or Dragovi\'{c} and Radnovi\'{c}~\cite{DR2011}, Chapter~5). Though
it is less elegant, it yields explicit formulas which allow us to prove
somewhat more than $a_* \not=b_*$.

\begin{lem}[]\label{lem.bary}
If $0<c<1$, the barycenters $a_*$ and $b_*$ of the two symmetric Poncelet
polygons $a_1,a_2,\ldots,a_n$ and $b_1,b_2,\ldots,b_n$ satisfy $0<a_*<b_*<1$.
\end{lem}

\noindent
\par\noindent
\begin{minipage}[T]{0.55\linewidth}
\emph{Proof.}
Denote by $p_k=\mathrm{e}^{2\mathrm{i}\varphi_k}$, $k=0,1,\ldots,n$, the vertices
of a Poncelet $n$-gon with $p_n=p_0$ and
\[
\varphi_0<\varphi_1< \ldots< \varphi_{n} = \varphi_0+\pi.
\]
By elementary geometry (see Figure~3, the blue angle is $\varphi_k-\varphi_{k-1}$
and the green angle is $\varphi_k+\varphi_{k-1}$), we get
\begin{equation} \label{eq.Trig}
\cos(\varphi_k-\varphi_{k-1}) - c\, \cos(\varphi_k+\varphi_{k-1}) = r,
\end{equation}
and hence
\[
(R-c)\,\cos\varphi_k\cos\varphi_{k-1} + (R+c)\, \sin\varphi_k\sin\varphi_{k-1} = r.
\]
\end{minipage}
\hfill
\begin{minipage}[T]{0.4\linewidth}
\begin{overpic}[width=\linewidth,tics=5]{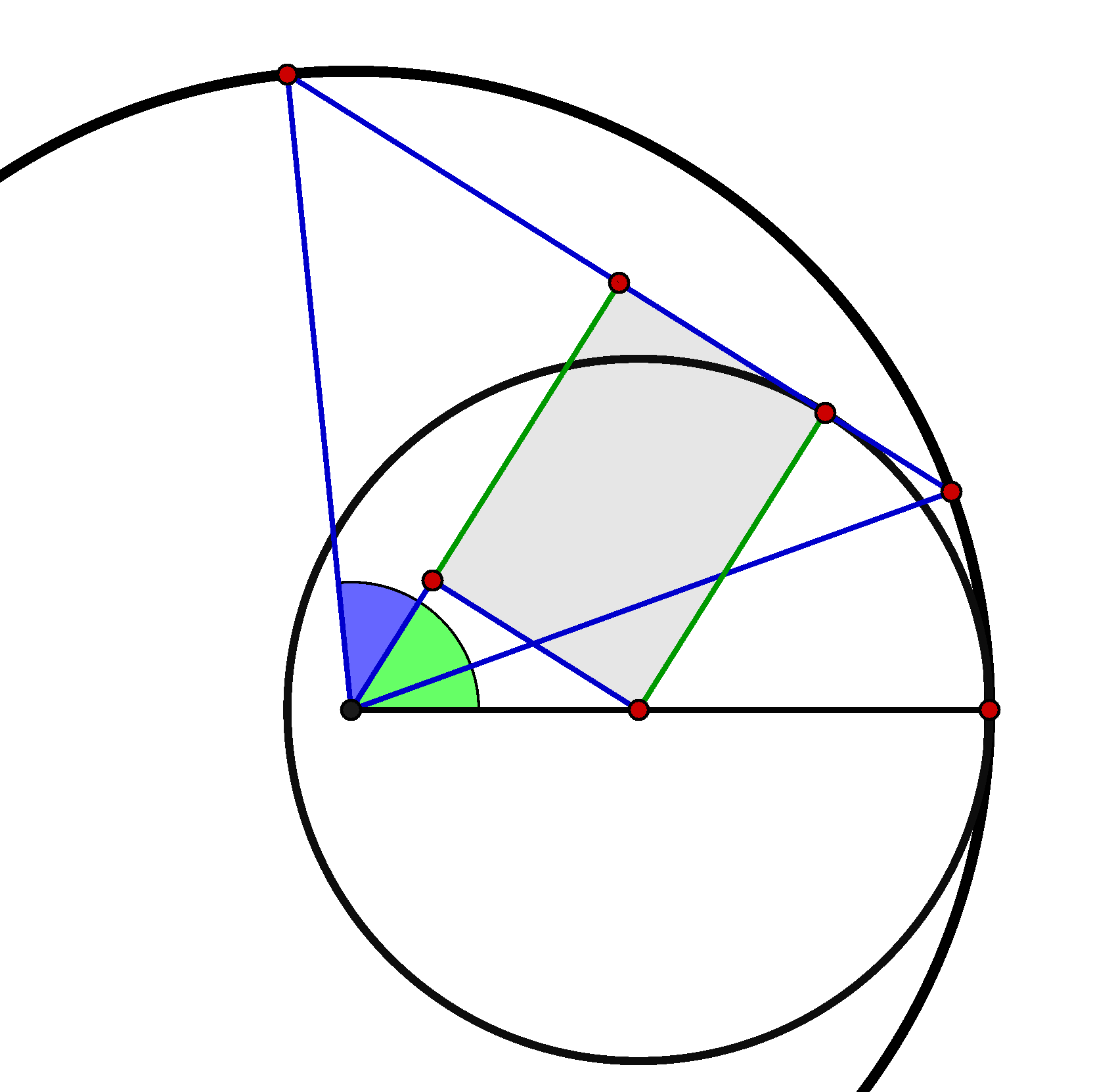}\\[1ex]
\fontsize{12pt}{14pt}\selectfont
\put(32,32){\makebox(0,0)[ct]{$0$}}
\put(58,32){\makebox(0,0)[ct]{$c$}}
\put(93,35){\makebox(0,0)[lc]{$1$}}
\put(90,55){\makebox(0,0)[lc]{$\mathrm{e}^{2\mathrm{i}\varphi_{k-1}}$}}
\put(27,95){\makebox(0,0)[cb]{$\mathrm{e}^{2\mathrm{i}\varphi_k}$}}
\put(49,58){\makebox(0,0)[lc]{$r$}}
\put(27,65){\makebox(0,0)[rc]{$1$}}
\end{overpic}
Figure~3: Derivation of equation~\eqref{eq.Trig}
\end{minipage}

\bigskip\par
\noindent
Setting $\psi_k:=\pi/2-\varphi_k$ we get%
\footnote{This substitution is needed to make the factor in front of the
sine functions less than one.}
\begin{equation} \label{eq.addthm1}
\cos\psi_k\cos\psi_{k-1} + \frac{1-c}{1+c}\, \sin\psi_k\sin\psi_{k-1} = \frac{r}{1+c}.
\end{equation}
In order to interpret this equation as an addition theorem for elliptic functions
we introduce the \emph{parameter}%
\footnote{Which is related to the more common \emph{elliptic modulus} $k$ by $m=k^2$.}
\begin{equation} \label{eq.DefM}
m := \frac{4c}{(1+c)^2-r^2},
\end{equation}
Defining
\[
t_k := \int_{0}^{\psi_k} \frac{d\theta}{\sqrt{1-m \sin^2\theta}}, \qquad k=0,\ldots,n,
\]
we have $\cos\psi_k = \cn t_k$ and $\sin\psi_k = \sn t_k$, with Jacobi's elliptic
functions {\it cosinus amplitudinis} and {\it sinus amplitudinis}, respectively.
This substitution converts \eqref{eq.addthm1} to
\begin{equation} \label{eq.addthm2}
\cn t_k\cn t_{k-1} + \frac{1-c}{1+c}\, \sn t_k \sn t_{k-1} = \frac{r}{1+c}.
\end{equation}
The functions $\sn$ and $\cn$ are doubly periodic; $\sn$ has fundamental periods
$4K$ and $2\mathrm{i}K'$, fundamental periods of $\cn$ are $4K$ and
$2K+2\mathrm{i}K'$, respectively. Here $K=K(m)$ is the \emph{complete elliptic
integral}
\[
K(m)= \int_{0}^{\pi/2}\frac{d\theta}{\sqrt{1-m \sin^2 \theta}}, \quad \text{and}
\quad K'=K(1-m).
\]
Figure~4 shows enhanced phase plots  of $\sn$ (left) and $\cn$ (right) with parameter
$m\approx 0.686$ (corresponding to $c=.1$ and $n=5$). The functions are depicted
in a domain somewhat larger than $-3K<\mathrm{Re}\,z<3K$, $-K'<\mathrm{Im}\,z<K'$.
The white lines are boundaries of a fundamental parallelogram.
\begin{center}
\includegraphics[width= 0.45\linewidth]{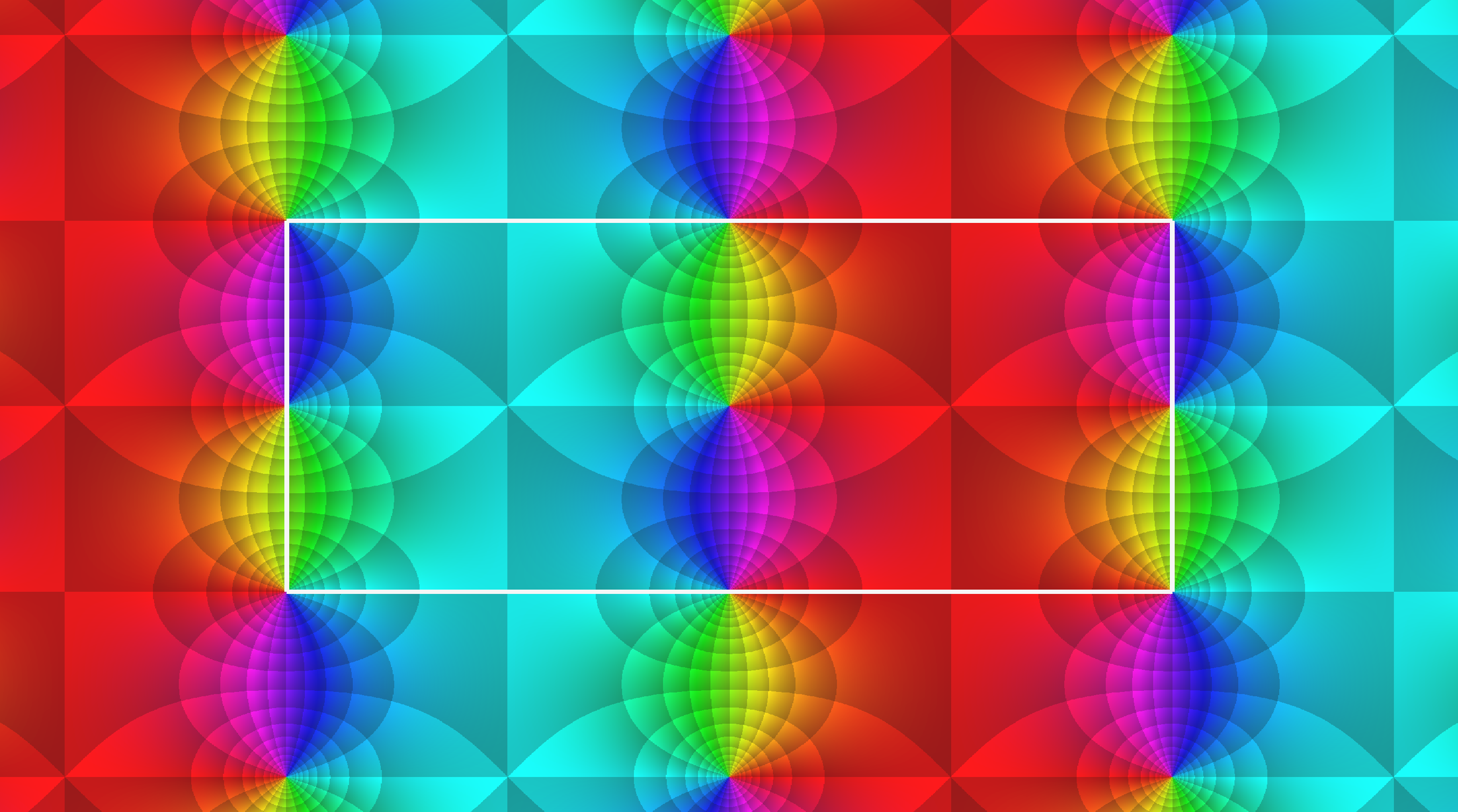}
\hfill
\includegraphics[width= 0.45\linewidth]{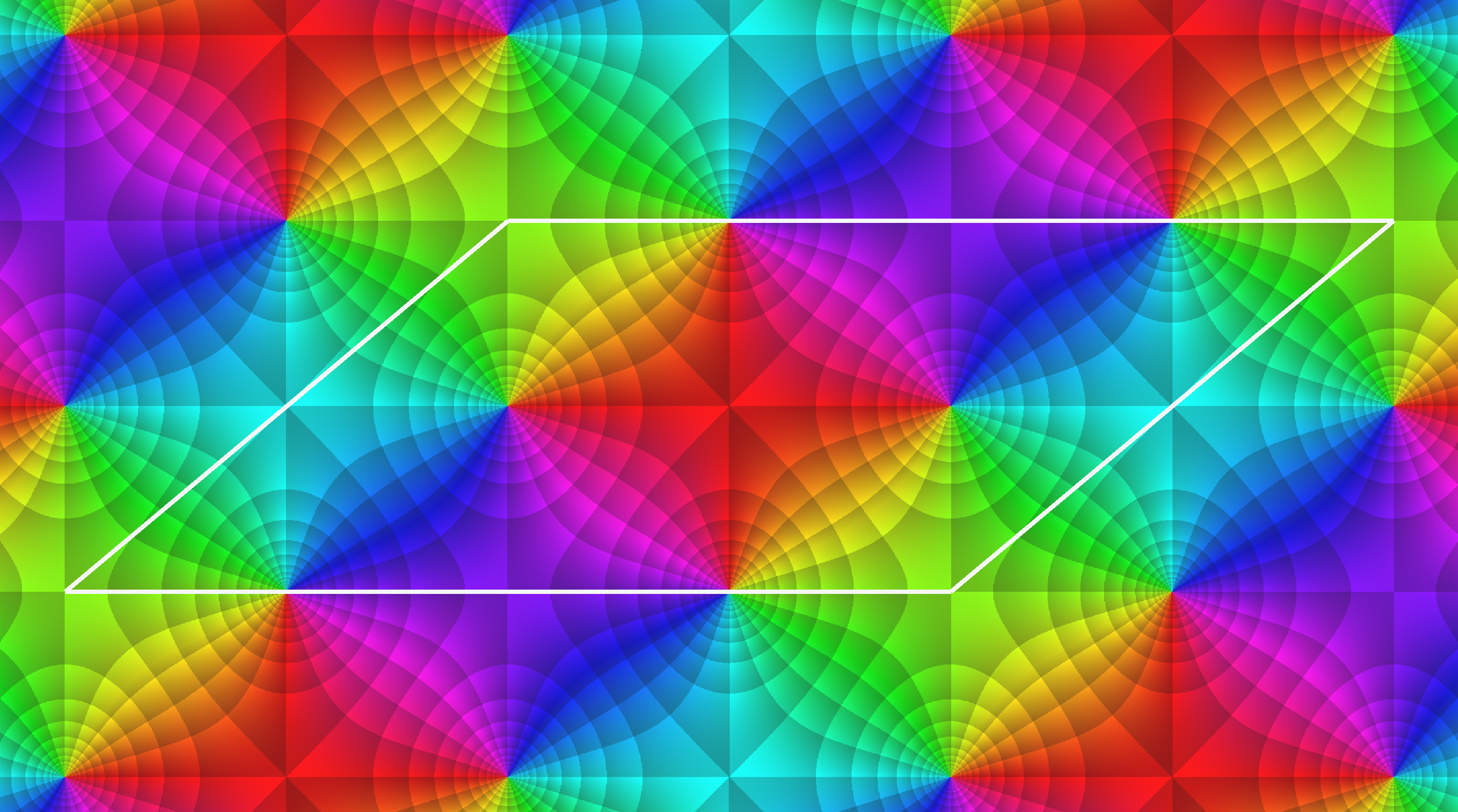}\\[1ex]
Figure~4: Phase plots of the functions $\sn$ (left) and $\cn$ (right)
\end{center}
Using the addition theorems for the functions $\cn$ and $\sn$, one can easily
verify the identity
\begin{equation} \label{eq.addthm3}
\cn (u+v)\cdot\cn u + \sqrt{1-m \sn ^2 v}\cdot\sn (u+v)\cdot\sn u = \cn v, \qquad
u,v \in \mathbb{C}.
\end{equation}
Let $s\in(0,K)$ be such that $\cn s = r/(1+c)$. Then the definition \eqref{eq.DefM}
of $m$ yields
\[
1-m\cdot\sn^2 s = 1-m\cdot(1-\cn^2 s)= \frac{(1-c)^2}{(1+c)^2}.
\]
After substituting $u:=t_{k-1}$ and $v:=s$ into \eqref{eq.addthm3} we arrive at
\begin{equation} \label{eq.addthm4}
\cn (t_{k-1}+s)\cdot\cn t_{k-1} + \frac{1-c}{1+c}\cdot\sn (t_{k-1}+s)
\cdot\sn t_{k-1} = \frac{r}{1+c}.
\end{equation}
In other words: \eqref{eq.addthm2} is satisfied for $t_{k}=t_{k-1}+s$.
Assuming for the moment that this indeed holds for all $k$, we get
\begin{equation} \label{eq.ExplicitTk}
t_k=t_0+k\,s,\qquad k=1,\ldots,n.
\end{equation}
Since the Poncelet polygon must be closed, $p_0=p_n$, we must have
$t_n=t_0+2\kappa\,K$ with some positive integer $\kappa$.
A little thought shows that $\kappa$ is the \emph{wrapping number} of the
polygon about the inner circle. In the case at hand we are interested
in solutions with $\kappa=1$, and hence
\begin{equation} \label{eq.CondD}
s=2K/n.
\end{equation}
Note that \eqref{eq.CondD} together with $\cn s=r/(1+c)$
implicitly determines the radius $r$ of the inner circle.
It can now be verified that \eqref{eq.ExplicitTk} is indeed the unique
solution (for fixed $t_0$) of \eqref{eq.addthm2} we are looking for.
Summarizing we get explicit formulas for the vertices,
\begin{equation} \label{eq.ExplicitPk}
p_k(t) = \mathrm{e}^{2\mathrm{i}\varphi_k}
= \mathrm{e}^{2\mathrm{i}(\pi/2-\psi_k)}
= - (\cos\psi_k-\mathrm{i}\sin\psi_k)^2
= - (\cn t_k - \mathrm{i}\,\sn t_k)^2, \quad t_k=t+k \,s.
\end{equation}
In what follows we consider the points $p_1,\ldots,p_n$ and their barycenter
$p_*$ as functions of the real parameter $t=t_0\in \mathbb{R}$.
Motivated by \eqref{eq.NRD2}, we are only interested in the real part of $p_*$,
\[
\varrho(t):=\mathrm{Re}\,p_*(t)=\frac{1}{n} \sum_{k=1}^{n} \mathrm{Re}\,p_k(t)
=\frac{1}{n} \sum_{k=1}^{n} \sigma(t+ks), \quad
\sigma(z):= 1-2\cn^2 z, \quad t \in \mathbb{R},\ z \in\mathbb{C}.
\]
The function $\sigma$ is an elliptic function with fundamental periods $2K$
and $2\mathrm{i}K'$. It has order $2$, with a double pole at $\mathrm{i}K'$.
A phase plot of $\sigma$ in the square $|\mathrm{Re}\,z|<K$, $|\mathrm{Im}\, z|<K'$
is depicted in Figure~5 on the left.

Since $\varrho$ is the sum of $n$ translations $t \mapsto t-2K\,k/n$
of $\sigma$, it extends to $\mathbb{C}$ as a meromorphic function with
periods $\omega_1=2K/n=s$ and $\omega_2=2\mathrm{i}K'$.
Let $\Lambda := \mathbb{Z}\,\omega_1 + \mathbb{Z}\,\omega_2$ denote the
period lattice of $\varrho$. Because the poles of $\sigma$ do not cancel in
the summation, $\varrho$ has (double) poles exactly at the points $\mathrm{i}K'+\Lambda$.
It follows that $\varrho$ is a non-constant (!) elliptic function of
order $2$. A phase plot of $\varrho$ in the square
$|\mathrm{Re}\,z|<K$, $|\mathrm{Im}\, z|<K'$ is shown in Figure~5, right.
The black and the white rectangles are boundaries of fundamental domains,
the black line is chosen such that it does not meet zeros or poles of
$\varrho$ and $\varrho'$.
\begin{center}
\includegraphics[width=0.45\linewidth]{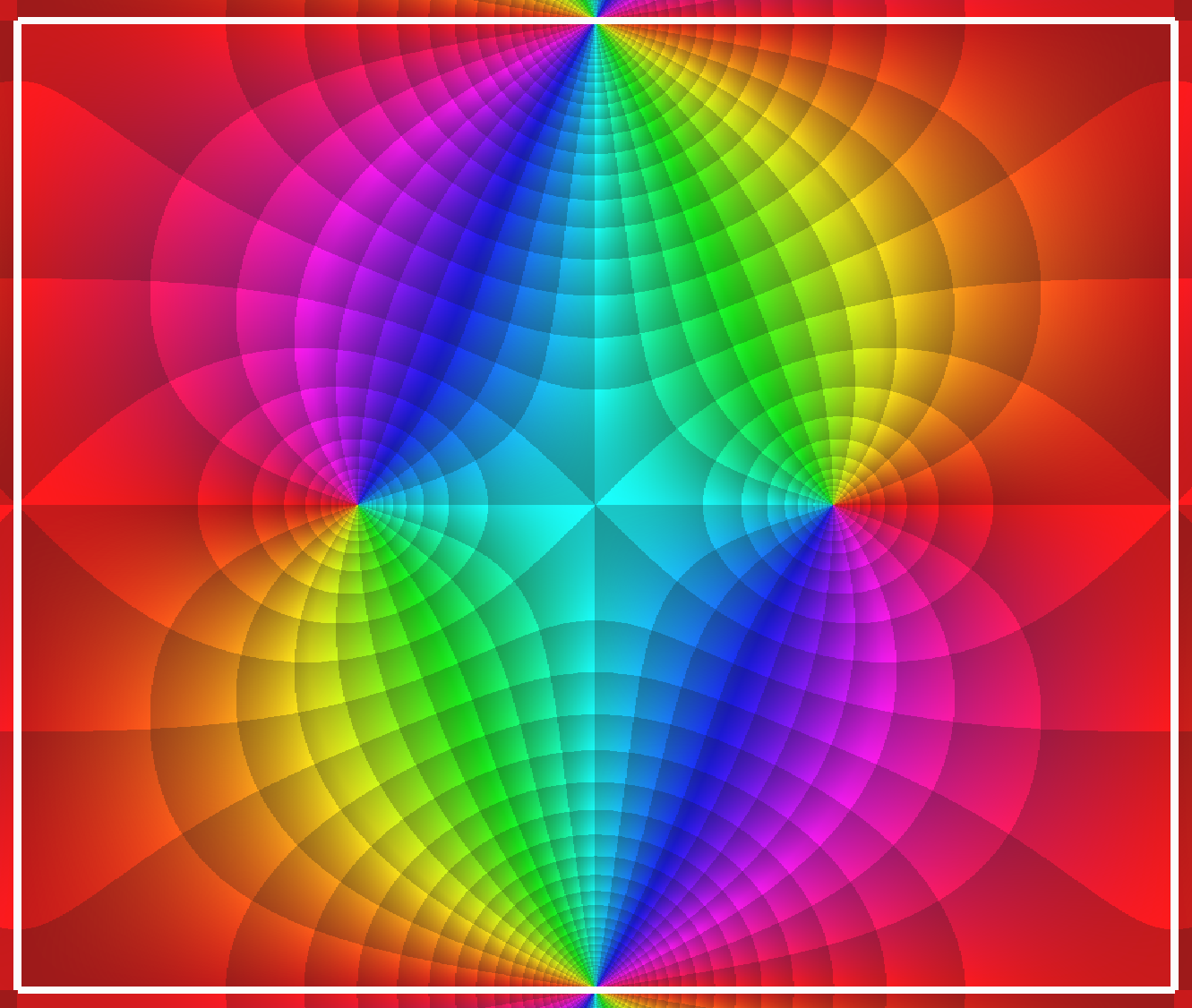}
\hfill
\includegraphics[width=0.45\linewidth]{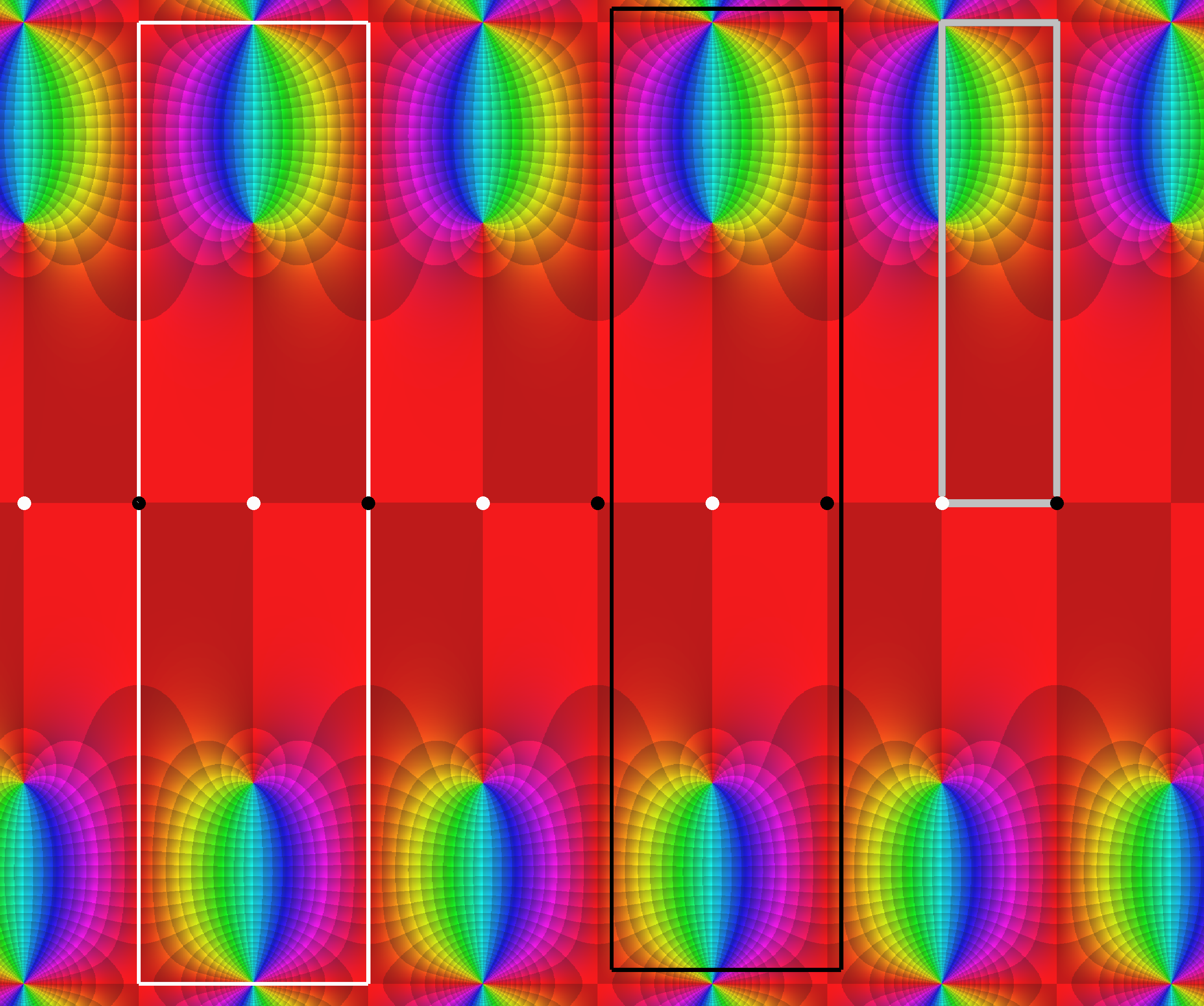}\\[1ex]
Figure~5: Phase plots of the functions $\sigma$ (left) and $\varrho$ (right)
\end{center}
Referring to the Schwartz-Tabachnikov result, we could finish the proof
here: If the locus of the barycenters were a point, the function $\sigma$
would be constant (on the real line, and thus in the entire plane),
which is not the case because it has poles.

For our more ambitious goal to prove the inequality $a_*<b_*$
we need some specific properties of $\varrho$.
Though these can be read off from the phase plot on the right-hand side
of Figure~5, we derive them from well-known facts about elliptic functions
and some symmetry arguments.

Let us first count the number of zeros and poles of $\varrho$ and
$\varrho'$ in a fundamental domain $\Omega$ of $\varrho$, chosen such
that none of these points lie on the boundary of $\Omega$
(the black line in the right image of Figure~5 bounds such a domain).
The function $\varrho$ has exactly one pole in $\Omega$, and its
multiplicity is two, so that $\varrho'$ has exactly one pole with multiplicity
three. By Liouville's Theorem $\varrho$ has two zeros, while $\varrho'$ has
three zeros (counted with multiplicity).

It is clear that $\varrho$ is real on the ``vertical'' lines
$\mathrm{Im}\,z=K' \,\mathbb{Z}$. Using symmetry properties, it can easily be
seen that it is also real on the ``horizontal'' lines
$\mathrm{Re}\,z=(s/2)\,\mathbb{Z}$. On those of these lines which do not
contain poles of $\varrho$ we find at least two zeros of $\varrho'$ in $\Omega$
(corresponding to maxima and minima of the real valued function $\varrho$),
while on those lines which contain double(!) poles of $\varrho$
there must be at least one zero of $\varrho'$ in $\Omega$. This
gives a total number of at least $2+1$ zeros on the vertical and another $2+1$
zeros on
the horizontal lines (in $\Omega$). Since $\varrho'$ has only $3$ zeros in $\Omega$,
these points must all be located at the crossings $0+\Lambda$,
$K/n+\Lambda$, $K/n+\mathrm{i}K' +\Lambda$ of the horizontal and the vertical
lines. Moreover, all zeros of $\varrho'$ must be simple and different from the
zeros of $\varrho$, i.e., they are \emph{saddle points} of $\varrho$.

Starting at the pole $\mathrm{i}K'$, we now walk along a rectangular path $R$
with vertices at $\mathrm{i}K'$, $0$, $K/n$ and $K/n+\mathrm{i}K'$ until we return to
$\mathrm{i}K'$ (the grey line in Figure~5 is a translation of $R$ by $4s$).
The function $\varrho$ is real on $R$, starts out from $-\infty$ at the
beginning, and must be strictly monotone as long as we do not meet zeros of
$\varrho'$. Taking into account the behavior of analytic functions at
their saddle points (see~\cite{Wgt2010},~\cite{Wgt2012}),
we get that $\varrho$ is strictly increasing along the whole path $R$.
Because we meet $0$ earlier than $K/n$ this implies
$a_*=\varrho(0)<\varrho(K/n)=b_*$.
\par\hfill$\Box$
\bigskip\par

\noindent
{\bf Acknowledgement.} We thank Gunter Semmler for bringing reference~\cite{SchT}
to our attention.

\bibliographystyle{amsplain}
\bibliography{master}

\providecommand{\bysame}{\leavevmode\hbox to3em{\hrulefill}\thinspace}
\providecommand{\MR}{\relax\ifhmode\unskip\space\fi MR }
\providecommand{\MRhref}[2]{%
  \href{http://www.ams.org/mathscinet-getitem?mr=#1}{#2}
}
\providecommand{\href}[2]{#2}
\begin{thebibliography}{10}

\bibitem{Bol2018}
V.~Bolotnikov, \emph{Boundary interpolation by finite {B}laschke products},
  Complex analysis and dynamical systems, Trends Math.,
  Birkh\"{a}user/Springer, Cham, 2018, pp.~39--65.

\bibitem{CaPh1965}
D.~G. Cantor and R.~R. Phelps, \emph{An elementary interpolation theorem},
  Proc. Amer. Math. Soc. \textbf{16} (1965), 523--525.

\bibitem{ChaGoPaRo18}
I.~Chalendar, P.~Gorkin, J.~R. Partington, and W.~T. Ross, \emph{Clark measures
  and a theorem of {R}itt}, Math. Scand. \textbf{122} (2018), no.~2, 277--298.

\bibitem{DGSV}
U.~Daepp, P.~Gorkin, A.~Shaffer, and K.~Voss, \emph{Finding ellipses}, Carus
  Mathematical Monographs, vol.~34, MAA Press, Providence, RI, 2018, What
  Blaschke products, Poncelet's theorem, and the numerical range know about
  each other.

\bibitem{DaGoVo}
U.~Daepp, P.~Gorkin, and K.~Voss, \emph{Poncelet's theorem, {S}endov's
  conjecture, and {B}laschke products}, J. Math. Anal. Appl. \textbf{365}
  (2010), no.~1, 93--102.

\bibitem{DR2011}
V.~Dragovi\'c and M.~Radnovi\'c, \emph{{Poncelet porisms and beyond. Integrable
  billiards, hyperelliptic Jacobians and pencils of quadrics}}, Basel:
  Birkh\"auser, 2011.

\bibitem{GaMaRo}
S.~R. Garcia, J.~Mashreghi, and Ross~W. T., \emph{Finite {B}laschke products
  and their connections}, Springer, Cham, 2018.

\bibitem{GWW16}
H.-L. Gau, K.-Z. Wang, and P.~Y. Wu, \emph{Circular numerical ranges of partial
  isometries}, Linear Multilinear Algebra \textbf{64} (2016), no.~1, 14--35.

\bibitem{GauWu98}
H.-L. Gau and P.~Y. Wu, \emph{Numerical range of {$S(\phi)$}}, Linear
  Multilinear Algebra \textbf{45} (1998), no.~1, 49--73.

\bibitem{Gla2008}
Ch. Glader, \emph{{Minimal degree rational unimodular interpolation on the unit
  circle}}, {ETNA, Electron. Trans. Numer. Anal.} \textbf{30} (2008), 88--106.

\bibitem{GoWa17}
P.~Gorkin and N.~Wagner, \emph{Ellipses and compositions of finite {B}laschke
  products}, jmaa \textbf{445} (2017), no.~2, 1354--1366.

\bibitem{Gri}
Ph.~A. Griffiths, \emph{Complex analysis and algebraic geometry}, Bull. Amer.
  Math. Soc. (N.S.) \textbf{1} (1979), no.~4, 595--626.

\bibitem{GusRa}
K.~E. Gustafson and D.~K.~M. Rao, \emph{Numerical range. {T}he field of values
  of linear operators and matrices}, Springer, New York, 1997.

\bibitem{HJ2}
R.~A. Horn and C.~R. Johnson, \emph{Topics in matrix analysis}, Cambridge
  University Press, Cambridge, 1994, Corrected reprint of the 1991 original.

\bibitem{Jac1828}
C.~G.~J. {Jacobi}, \emph{{Ueber die Anwendung der elliptischen Transcendenten
  auf ein bekanntes Problem der Elementargeometrie.}}, {J. Reine Angew. Math.}
  \textbf{3} (1828), 376--389 (in German).

\bibitem{JoRu1987}
W.~B. Jones and S.~Ruscheweyh, \emph{Blaschke product interpolation and its
  application to the design of digital filters}, Constr. Approx. \textbf{3}
  (1987), no.~4, 405--409.

\bibitem{Kipp1951}
R.~Kippenhahn, \emph{{\"Uber den Wertevorrat einer Matrix}}, {Math. Nachr.}
  \textbf{6} (1951), 193--228 (in German).

\bibitem{Ki08}
\bysame, \emph{On the numerical range of a matrix}, Linear Multilinear Algebra
  \textbf{56} (2008), no.~1-2, 185--225, Translated from the German by Paul F.
  Zachlin and Michiel E. Hochstenbach.

\bibitem{Mir}
B.~Mirman, \emph{Numerical ranges and {P}oncelet curves}, Linear Algebra Appl.
  \textbf{281} (1998), 59--85.

\bibitem{SchT}
R.~Schwartz and S.~Tabachnikov, \emph{Centers of mass of {P}oncelet polygons,
  200 years after}, Math. Intelligencer \textbf{38} (2016), no.~2, 29--34.

\bibitem{SeWe2006}
G.~Semmler and E.~Wegert, \emph{Boundary interpolation with {B}laschke products
  of minimal degree}, Comput. Methods Funct. Theory \textbf{6} (2006), no.~2,
  493--511.

\bibitem{SSW}
I.~M. Spitkovsky, I.~Suleiman, and E.~Wegert, \emph{The {G}au-{W}ang-{W}u
  conjecture on partial isometries holds in the 5-by-5 case}, arXiv:2108.07024
  [math.FA] (2021), 1--10.

\bibitem{Wgt2010}
E.~Wegert, \emph{{Phase diagrams of meromorphic functions}}, {Comput. Methods
  Funct. Theory} \textbf{10} (2010), no.~2, 639--661 (English).

\bibitem{Wgt2012}
\bysame, \emph{{Visual complex functions. An introduction with phase
  portraits}}, Basel: Birkh\"auser, 2012.

\bibitem{WeSe}
E.~Wegert and G.~Semmler, \emph{Phase plots of complex functions: a journey in
  illustration}, Notices Amer. Math. Soc. \textbf{58} (2011), no.~6, 768--780.

\end{thebibliography}

\end{document}